\documentclass[11pt]{article} 

 \usepackage{booktabs}       
 \usepackage{amsfonts}       
 \usepackage{nicefrac}       
 \usepackage{microtype}      
 
\usepackage{booktabs} 

\usepackage[toc,page]{appendix}

\usepackage{enumerate,wrapfig,placeins}
\usepackage{bbm}
\usepackage{graphicx}
\usepackage{enumitem}  
\usepackage{amsmath,amsthm,amsfonts,amssymb}  
\usepackage{pifont}

\usepackage{geometry}

\pdfoutput=1  

\usepackage{hyperref}

\usepackage{amsmath, amssymb, graphicx, url}

\usepackage{optidef}
\usepackage{graphics} 
\graphicspath{{figures/}}

\usepackage{url}            

\usepackage{wrapfig}

\usepackage{caption}
\usepackage{marginnote}
\usepackage{float} 

\usepackage{textgreek,upgreek,bm}
\usepackage{ulem} 

     \makeatletter
\newcommand\gobblepars{%
    \@ifnextchar\par%
        {\expandafter\gobblepars\@gobble}%
        {}}
\makeatother

\def\mindex#1{\index{#1}}

\def\Obj{L}

\def\Tdiff{\mathcal{D}}

\def\fee{\upphi}

\def\uH{\underline{H}}
\def\uQ{\underline{Q}}

\newcommand{\bbblot}{\raise1pt\hbox{\vrule height .4ex width .4ex depth .05ex}}

\long\def\defbox#1{\framebox[.9\hsize][c]{\parbox{.85\hsize}{%
\parindent=0pt
\baselineskip=12pt plus .1pt      
\parskip=6pt plus 1.5pt minus 1pt 
 #1}}}

\long\def\beginbox#1\endbox{\subsection*{}%
\hbox{\hspace{.05\hsize}\defbox{\medskip#1\bigskip}}%
\subsection*{}}

\def\endbox{}

 \def\archival#1{} 

\def\FRAC#1#2#3{\genfrac{}{}{}{#1}{#2}{#3}}

\def\ddtp{{\mathchoice{\FRAC{1}{d^{\hbox to 2pt{\rm\tiny +\hss}}}{dt}}%
{\FRAC{1}{d^{\hbox to 2pt{\rm\tiny +\hss}}}{dt}}%
{\FRAC{3}{d^{\hbox to 2pt{\rm\tiny +\hss}}}{dt}}%
{\FRAC{3}{d^{\hbox to 2pt{\rm\tiny +\hss}}}{dt}}}}

\def\ddyp{{\mathchoice{\FRAC{1}{d^{\hbox to 2pt{\rm\tiny +\hss}}}{dy}}%
{\FRAC{1}{d^{\hbox to 2pt{\rm\tiny +\hss}}}{dy}}%
{\FRAC{3}{d^{\hbox to 2pt{\rm\tiny +\hss}}}{dy}}%
{\FRAC{3}{d^{\hbox to 2pt{\rm\tiny +\hss}}}{dy}}}}

\def\half{{\mathchoice{\FRAC{1}{1}{2}}%
{\FRAC{1}{1}{2}}%
{\FRAC{3}{1}{2}}%
{\FRAC{3}{1}{2}}}}

\def\darrow{\buildrel{\rm dist}\over\longrightarrow}

\newsavebox{\junk}
\savebox{\junk}[1.6mm]{\hbox{$|\!|\!|$}}

\def\argmin{\mathop{\rm arg{\,}min}}
\def\argmax{\mathop{\rm arg{\,}max}}

\def\state{{\sf X}}

\def\ustate{{\sf U}} 
\def\ystate{{\sf Y}} 
\def\zstate{{\sf Z}}

\def\ystate{{\sf Y}}

\def\bfmath#1{{\mathchoice{\mbox{\boldmath$#1$}}%
{\mbox{\boldmath$#1$}}%
{\mbox{\boldmath$\scriptstyle#1$}}%
{\mbox{\boldmath$\scriptscriptstyle#1$}}}}

\def\bfmU{\bfmath{U}}

\def\bfmX{\bfmath{X}}

\def\bfmY{\bfmath{Y}}

\def\bfmhhaY{\bfmath{\hhaY}} 
\def\bfmhhaY{\hbox to 0pt{$\widehat{\bfmY}$\hss}\widehat{\phantom{\raise 1.25pt\hbox{$\bfmY$}}}}

\def\bfmZ{\bfmath{Z}}

\def\hatheta{{\hat\theta}}

\def\haX{\widehat X}

\def\tilA{\tilde{A}}

\def\tiltheta{{\tilde \theta}}

\def\tilg{\tilde g}

\def\clE{{\cal E}}
\def\clF{{\cal F}}

\def\clI{{\cal I}}

\def\clR{{\cal R}}
\def\clS{{\cal S}}

\def\clI{{\cal I}}

\def\eqdef{\mathbin{:=}}

\def\Prob{{\sf P}}

\def\Expect{{\sf E}}

\def\lgmath#1{{\mathchoice{\mbox{\large #1}}%
{\mbox{\large #1}}%
{\mbox{\tiny #1}}%
{\mbox{\tiny #1}}}}

\def\One{{\mathchoice{\lgmath{\sf 1}}%
{\mbox{\sf 1}}%
{\mbox{\tiny \sf 1}}%
{\mbox{\tiny \sf 1}}}}

\def\ind{\bbbone}
  
 \def\epsy{\varepsilon}

\def\varble{\,\cdot\,}

\def\formtmp#1#2{{\vskip12pt\noindent\fboxsep=0pt\colorbox{#1}{\vbox{\vskip3pt\hbox to \textwidth{\hskip3pt\vbox{\raggedright\noindent\textbf{#2\vphantom{Qy}}}\hfill}\vspace*{3pt}}}\par\vskip2pt%
\noindent\kern0pt}}

\newenvironment{programcode}[1]{\ignorespaces\def\stmtopen##1{##1}%
\pagebreak[3]%
\formtmp{programcode}{#1}%
\endlinechar=-1\relax%
\nopagebreak[4]}{%
\noindent\textcolor{programcode}{\rule{\columnwidth}{1pt}}\vskip1pt\par\addvspace{\baselineskip}%
\endlinechar=13}
 
\def\barg{{\overline {g}}}

\def\barr{{\overline {r}}}

\def\barA{{\bar{A}}}

\def\bartheta{{\overline{\theta}}}

\def\barbeta{{\bar{\beta}}}

{\end{list}}

\def\ass(#1:#2){(#1\ref{#1:#2})}

\def\ritem#1{
\item[{\sf \ass(\current_model:#1)}]
}

\newenvironment{recall-ass}[1]{%
\begin{description}
\def\current_model{#1}}{
\end{description}
}

\def\sq{\hbox{\rlap{$\sqcap$}$\sqcup$}}
\def\qed{\ifmmode\sq\else{\unskip\nobreak\hfil
\penalty50\hskip1em\null\nobreak\hfil\sq
\parfillskip=0pt\finalhyphendemerits=0\endgraf}\fi}

\newcommand{\blot}{\vrule height 1.1ex width .9ex depth -.1ex }
\def\qedb{\ifmmode\blot\else{\vspace{-.2cm}\unskip\nobreak\hfil
\penalty50\hskip1em\null\nobreak\hfil\blot
\parfillskip=0pt\finalhyphendemerits=0\endgraf}\fi}

\newcounter{rmnum}
\newenvironment{romannum}{\begin{list}{{\upshape (\roman{rmnum})}}{\usecounter{rmnum}
\setlength{\leftmargin}{8pt}
\setlength{\rightmargin}{6pt}
\setlength{\itemindent}{0pt}
}}{\end{list}}

\newcounter{anum}

\newcommand{\field}[1]{\mathbb{#1}}

\def\Re{\field{R}}

\def\Prob{{\sf P}}

\def\Expect{{\sf E}}

\def\transpose{{\intercal}}

\def\st{\text{\rm s.t.\,}}

\def\argmin{\mathop{\rm arg\, min}}
\def\ind{\hbox{\large \bf 1}}

\def\epsy{\varepsilon}
\def\varble{\,\cdot\,}

\def\haJ{\widehat J}

\def\haY{\widehat{Y}}

\def\hhaY{\hbox to 0pt{$\haY$\hss}\widehat{\phantom{\raise 1.25pt\hbox{Y}}}}

\def\haY{\widehat Y}

\def\wham#1{\smallbreak\pagebreak[3]%
	\noindent\textbf{#1}\ \ \gobblepars}
	
	\def\whamit#1{\smallbreak\pagebreak[3]%
	\noindent\textit{#1}\ \ \gobblepars}

\graphicspath{{figures/}}

\usepackage{color}

\newlength{\noteWidth}
\long\def\notes#1{\ifinner
	{\tiny #1}
	\else
	\marginpar{\parbox[t]{\noteWidth}{\raggedright\tiny #1}}
	\fi}

\def\notes#1{\typeout{See notes!}}

\def\sfb#1{}

\def\choreL#1{\notes{\color{green} FL:  #1}}

\newtheorem{theorem}{Theorem}[section]

\newtheorem{proposition}[theorem]{Proposition}
\newtheorem{lemma}[theorem]{Lemma}

\usepackage{cleveref}

\Crefname{corollary}{Corollary}{Corollaries}
\Crefname{eqnarray}{eq.}{eqs.}
\Crefname{equation}{eq.}{eqs.}

\Crefname{figure}{Fig.}{Figs.}
\Crefname{tabular}{Tab.}{Tabs.}
\Crefname{table}{Tab.}{Tabs.}
\Crefname{proposition}{Prop.}{Propositions}
\Crefname{theorem}{Thm.}{Thms.}
\Crefname{definition}{Def.}{Defs.} 
\Crefname{section}{Section}{Sections}
\Crefname{lemma}{Lemma}{Lemmas}
\Crefname{assumption}{Assumption}{Assumptions}

\def\Obj{\Upgamma}  
\def\barObj{\bar{\Obj}}

 \def\SigmaTheta{\Sigma_\uptheta}

\def\thetaRegion{\Uptheta}   

\def\psisub#1{\psi_{(#1)}}
\def\upsisub#1{\underline{\psi}_{(#1)}}

\def\csub#1{c_{(#1)}}

\def\TdiffR{\widehat{\Tdiff}}

\def\objv{v} 

\def\constD{{d_{+}}} 

\def\tilA{\widetilde{A}}
\def\barW{\overline{W}}
\def\ystate{\textsf{Y}}
\def\barN{\overline{N}}

\def\tilbeta{\widetilde\beta}

\makeatletter
\newcommand{\ostar}{\mathbin{\mathpalette\make@circled\star}}
\newcommand{\make@circled}[2]{%
	\ooalign{$\m@th#1\smallbigcirc{#1}$\cr\hidewidth$\m@th#1#2$\hidewidth\cr}%
}
\newcommand{\smallbigcirc}[1]{%
	\vcenter{\hbox{\scalebox{0.77778}{$\m@th#1\bigcirc$}}}%
}
\makeatother

\title{ 	  Convex Q Learning in a Stochastic Environment:  Extended Version
}

\author{Fan~Lu%
	\thanks{F.~Lu is with the Department of Applied Mathmatics and Statistics, University of California, Santa Cruz, CA 95064, USA (e-mail: flu16@ucsc.edu).}
	\
 and  Sean~Meyn
	\thanks{S.~Meyn is with the Department of Electrical and Computer Engineering, University of Florida, Gainesville, FL 32611, USA (e-mail: meyn@ece.ufl.edu).
 Support from  ARO award  W911NF2010055 and   National Science Foundation award    
		   EPCN 1935389 is gratefully acknowledged.  }    
}

\begin{document}
\maketitle

\begin{abstract}

The paper introduces the first formulation of convex Q-learning for Markov decision processes with function approximation.   The algorithms and theory rest on a relaxation of a dual of Manne's celebrated linear programming characterization of optimal control.
The main contributions firstly concern properties of the relaxation, described as a deterministic convex program:  we identify conditions for a bounded solution,  and a significant relationship between the solution to the new convex program, and the solution to standard Q-learning.   The second set of contributions concern algorithm design and analysis:    
(i) A direct model-free method for approximating the convex program for Q-learning shares properties with its ideal.   In particular, a bounded solution is ensured subject to a simple property of the basis functions;  (ii)  The proposed algorithms are convergent and new techniques are introduced to obtain the rate of convergence in a mean-square sense;   
(iii) The approach can be generalized to a range of performance criteria, and it is found that variance can be reduced by considering ``relative'' dynamic programming equations;  
(iv) The theory is illustrated with an application to a classical inventory control problem.

\noindent
This is an extended version of an article to appear in the forthcoming IEEE Conference on Decision and Control.

\end{abstract}

\section{Introduction}
\label{qlearning}


The Q-learning algorithm introduced in \cite{wat89} is a highly celebrated approach to reinforcement learning, that has evolved over the past decades to form part of the solution to complex optimal control problems.    It was originally designed to compute the state-action value function (known as the Q-function).   This early work considered the discounted-cost optimal control problem for Markov decision processes (MDPs),    in the  \textit{tabular} setting  so that the function class spans all functions.     

The ultimate goal then and now is to approximate the Q-function within a restricted function class, notably neural networks, though much of the theory is restricted to a linearly parameterized function class.    Counterexamples show that conditions on the function class are required in general,  even in a linear function approximation setting  \cite{bai95,tsivan96,gor00}.
Criteria for stability based on sufficient exploration are contained in the recent work \cite{mey23}.


Moreover, when convergent, the limit of Q-learning or DQN solves a ``projected Bellman equation'' (see \eqref{e:projBE}), but we know little about the implication of this conclusion.   These concerns have motivated new ways of thinking about how to approximate a Q-function.

One alternative is GQ learning, based on a stochastic gradient descent algorithm with an objective similar to the mean-square Bellman error \cite{sutszemae08}.   Recently a more exact stochastic gradient descent algorithm was introduced in \cite{avrbordolpat21} with full stability analysis.      These results present a significant advancement but come with two drawbacks: the objective is non-convex, so there is no reason to believe the algorithm will converge to the global minimum.   Moreover, it remains difficult to interpret the solution of the global minimum.   If the Bellman error is small in an $L_2$ sense, where the $L_2$ norm depends on training data, what does this tell us about the performance of the ultimate feedback policy?

The linear programming (LP)  approach to optimal control pioneered by Manne  \cite{man60a} has inspired alternative approaches to RL and approximate dynamic programming.   
The earliest such work was found in   \cite{schsei85},  with error bounds appearing in
 \cite{farroy06,defvan04,lakbhasze17a}. 
Model-free algorithms appeared in \cite{mehmey09a,mehmeyneulu21,lumehmeyneu22} and \cite[Ch.~5]{CSRL}, where the term \textit{convex Q-learning} (CvxQ) was coined.
In parallel came \textit{logistic Q-learning} \cite{bascurkraneu21},  which solves a regularized dual of the LP in \cite{mehmeyneulu21}.     
There is however a gap in the settings:  CvxQ was developed for deterministic control systems,  while logistic Q-learning treats MDPs.   Also, the stochastic setting is so far restricted to tabular  \cite{bascurkraneu21} or  linearly factorable MDPs \cite{neuoko22}.      
Theory for CvxQ has few restrictions, beyond the limitation to deterministic control systems.

LP approaches are attractive because we obtain by design a convergent algorithm.  Moreover, the $L_\infty$-framework is more likely to lead to an interpretable solution,  since performance bounds on the resulting feedback policy can be obtained through Lyapunov function techniques  \cite{farroy06,lakbhasze17a}.  The main contributions are summarized here:

\wham{i)} Convex Q-learning for optimal control is introduced in a stochastic environment for the first time.   It is found that the constraint region is bounded subject to a persistence of excitation, generalizing the conclusions obtained recently for deterministic optimal control problems \cite{lumehmeyneu22}.
 Several approaches to approximating the solution to the convex program are proposed and analyzed.

\wham{ii)}   \Cref{t:cvxQ_Galerkin} implies a surprising connection between CvxQ and standard Q-learning.

\wham{iii)}   Techniques are introduced to obtain the rate of convergence in a mean-square sense---see \Cref{t:randLP}.

\wham{Comparison with existing literature.}    The new algorithms and some of the analysis might be anticipated from the theory for deterministic control systems in  \cite{lumehmeyneu22}. 
  \Cref{t:cvxQ_Galerkin} is new (and was initially surprising to us) even in the deterministic setting.
    The variance analysis surveyed in \Cref{t:randLP} is novel, resting on recent CLT theory from   \cite{ste01a} to obtain an exact formula for the asymptotic covariance.    Complementary results appeared in  \cite{petzil09}, motivated by MDP LP relaxations.   Conclusions in this prior work is based on i.i.d.\   samples of trajectories,  designed to permit application of Hoeffding's inequality to obtain sample complexity bounds for constraint-sampled LPs.      The   covariance formula in  \Cref{t:randLP} is similar to what is anticipated from  stochastic approximation (SA) theory, even though CvxQ falls outside of standard SA  recursions  ---see discussion following the proposition.

\wham{Organization:}  \Cref{s:qcp}  conditions foundations for the CvxQ algorithms introduced in \Cref{s:algs}.   
Theory for convergence rates is contained in \Cref{s:vartheory}.  The theory is illustrated in \Cref{s:numerics} with an application to a classical inventory control problem.

\bigskip

This paper is an extended version of \cite{lumey23a},    to appear in the 2023 IEEE Conference on Decision and Control.

%
%


\section{Q-learning Convex Programs}
\label{s:qcp}


The control model for which algorithm design and analysis is based on the standard  Markov Decision Process (MDP) 
with finite state space $\state$, finite input space $\ustate$,  and non-negative cost function $c : \zstate \to \Re_+$,  
with $\zstate \eqdef \state \times \ustate$.   

The state process is denoted $\bfmX =\{X(k):k\ge 0\}$, the
 input (or action) sequence   $\bfmU =\{U(k):k\ge 0\}$,  and the pair process $\bfmZ = (\bfmX,\bfmU)$.
The controlled transition matrix is denoted   $P_u$ for $u \in \ustate$, so that $\Prob\{X(k+1) = x' \mid X(k)=x, U(k)=u \}=P_u(x,x')$;  it 
  acts on functions $V:\state \to \Re$ via
$	P_u V(x) \eqdef \sum_{x^\prime \in \state} P_u(x, x^\prime) V(x^\prime)$.

\whamit{Choice of training input.}  The formulation of convex programs as well as theory surrounding algorithms is couched in a stationary setting:
 the input used for training is defined by a randomized stationary policy, so that  the joint process $\bfmZ$  is a time homogeneous Markov chain on   $\zstate$.
 It is assumed that $\bfmZ$ is uni-chain, with unique invariant pmf  denoted $\varpi$.   
 
 However, we often obtain faster training when using an   epsilon-greedy policy;  
 analysis of this more complex setting is beyond the scope of this paper.  

Theory is also restricted to the discounted-cost optimality criterion, with discount factor $\gamma \in (0, 1)$.  
The Q-function is the state-action value function,
\begin{equation*}
	\begin{aligned}
		Q^*(z) \eqdef \min_{\bfmU} \sum_{k=0}^{\infty} \gamma^k \Expect[c(Z(k)) | Z(0) = z] \,, \ \ z = (x, u) \in \zstate
	\end{aligned}
	\label{e:optValueFunc}
\end{equation*}
where the minimum is over all adapted input sequences.    It is   the unique solution to the Bellman equation,   
\begin{align}
	Q^*(z) = c(z) + \gamma P_u\uQ^*(x)
	\label{e:BellmanEq}
\end{align}
with $\uQ(x) \eqdef \min_u Q(x, u)$.
The optimal input is state feedback $U_k = \fee^*(X_k)$, using the ``$Q^*$-greedy" policy,
\begin{align}
	\fee^*(x) \in \argmin_{u\in \ustate} Q^* (x, u) \,, \qquad x \in \state
	\label{e:QgreedyPolicy}
\end{align}

\subsection{Convex programs for approximation}

Convex Q-learning algorithm is motivated by the classical LP characterization of optimal control problem due to Manne  \cite{man60a}.   The following is a simple corollary:


\begin{proposition}
\label[proposition]{t:lpcharact}
For any pmf $\mu$ on $\zstate$, the Q-function $Q^*$ is a solution to the  convex program,
	\begin{equation}
		\label{e:cvxQInit}
		\begin{aligned}
			\max_Q &\ \ \langle \mu, Q \rangle
			\\
			\st &\ \ Q(z) \le c(z) + \gamma  P_u \uQ(x)\, ,  \ \ z = (x, u)\in \zstate.
		\end{aligned}
	\end{equation}

\end{proposition}

\smallskip

The proof follows from verification that $Q \le Q^*$ whenever $Q$ is feasible,  which is  a standard Lyapunov function argument.  

This section is devoted to relaxations of 	\eqref{e:cvxQInit} that are model-based.   The conclusions motivate model-free algorithms and analysis in \Cref{s:algs}.



To obtain a convex program we restrict to a linear family:    $\{ Q^\theta(x, u) = \theta^\transpose \psi(x, u) : \theta \in \Re^d \}$, 
with $\psi: \state \times \ustate \to \Re^d$ the vector of basis functions, and based on an appropriate approximation obtain a policy in analog with \eqref{e:QgreedyPolicy}:
 \begin{align}
	\fee^\theta(x) \in \argmin_{u\in \ustate} Q^\theta(x, u)
	\label{e:QthetaPolicy}
\end{align}

The following is suggested by \eqref{e:cvxQInit},
\begin{equation}
	\max_\theta \ \ \langle \mu, Q^\theta \rangle
	\qquad
	\st \ \ Q^\theta(z) \le c(z) +\gamma P_u \uQ^\theta(x) 
 	\label{e:cvxQTheta}
\end{equation}
This is not practical in typical applications because 
it requires knowledge of the model,  
and there are so many constraints:  one for each $z=(x,u)\in\zstate$.

\wham{Galerkin relaxations}
Practical algorithms are obtained by expressing the constraints of  	\eqref{e:cvxQTheta}
 in  sample path form:   	
if a vector $\theta\in\Re^d$ is feasible, then the following inequality is valid for \textit{any} adapted input sequence:
with $\clF_{k} = \sigma(Z(i): i \le k)$ the filtration generated by the observations,
\begin{equation}
\begin{aligned}
&	\Expect\bigl[\Tdiff_{k + 1}(\theta) | \clF_k \bigr]  \ge 0 \,,   \quad \textit{for all $ k\ge 0$,} 
\\
&	\Tdiff_{k + 1}(\theta) \eqdef -Q^\theta(Z(k)) + c(Z(k)) + \gamma \uQ^\theta(X(k+1)) 
\end{aligned}
\label{e:TD}
\end{equation}

\notes{remove for CDC:}
\begin{lemma}
	\label[lemma]{t:CEpos}
	The following are equivalent for an integrable random variable $X$  and $\sigma$-field $\clF$:    
	\wham{\rm (i)}   $\haX \eqdef \Expect[X\mid \clF] \ge 0$ with probability one.
	
	\wham{\rm (ii)}   $\Expect[X \zeta] \ge 0$   for any $\clF$-measurable random variable $\zeta$ that is non-negative and bounded  (i.e., $\zeta\in L_\infty$).  
\end{lemma}

Armed with this lemma, 
a relaxation of \eqref{e:cvxQTheta} is obtained by specifying   a sequence of \textit{non-negative}  $\constD$-dimensional random vectors $\{\zeta_k:  k\ge 0\}$, with $\constD>d$.
Denote  $\barg(\theta) \eqdef \Expect_\varpi\bigl[  -\Tdiff_{k + 1}(\theta)\zeta_k \bigr]$.   A relaxation of  \eqref{e:cvxQTheta}
is then defined by
\begin{equation}
	\max_{\theta}
	\ \  \langle \mu, Q^\theta \rangle  \qquad
	\st  \ \
	\barg(\theta) \le  0.
	\label{e:ConvexQinfinity}
\end{equation}

An equivalent LP formulation is required for analysis.   
Let $\Upphi$ denote the set of all deterministic policies,   and for each $\fee\in \Upphi$ and $k\ge 0$ denote
\begin{equation}
\Tdiff_{k + 1}(\theta,\fee) =  \csub{k}  + \theta^\transpose \{ -\psisub{k}
+ \gamma  \psisub{k+1}^\fee  \}\, , 
\label{e:Tdiff_theta+fee}
\end{equation}
in which the following conventions will be used to save space when necessary:
$\csub{k} \eqdef c(Z(k)) $, and
\begin{equation*}
\psisub{k} \eqdef  \psi(X(k) , U(k ) ) \,, \ 
\psisub{k}^\fee \eqdef  \psi(X(k) , \fee(X(k)) ) 
\end{equation*}
Similar to \eqref{e:ConvexQinfinity}, we denote $\barg(\theta, \fee) \eqdef \Expect_\varpi\bigl[ -\Tdiff_{k + 1}(\theta, \fee)\zeta_k\bigr]$ for each $\fee \in \Upphi$.
\begin{proposition}
\label[proposition]{t:QCPLP}
Any solution to the \textit{Q-learning convex program} \eqref{e:ConvexQinfinity} is also a solution to the linear program, 
\begin{equation}
\label{e:ConvexQinfinityLP}
		\max_\theta
		\  \  \langle \mu, Q^\theta \rangle 
		\qquad
		\st   \ \ 
		\barg(\theta, \fee)  \le  0 \,, \ \   \  \fee\in\Upphi\,, 
\end{equation}
with identical optimal values.
\end{proposition}

\wham{Proof.}
For any $\fee \in \upphi$, we have
\begin{align*}
	\Expect_\varpi\bigl[ \Tdiff_{k + 1}(\theta,\fee) \zeta^i_k \bigr]  \ge \Expect_\varpi\bigl[ \Tdiff_{k + 1}(\theta) \zeta^i_k \bigr] \,, \ \ 1 \le  i \le \constD
\end{align*}
The proof is completed on recognizing that this lower bound   is achieved with $\fee = \fee^\theta$.
 \qed

Let $\thetaRegion  = \{ \theta \in \Re^d : \barg(\theta) \le 0 \}$ denote the constraint set for  \eqref{e:ConvexQinfinity}.
It is always non-empty since it contains the origin.     \Cref{t:constGeo} tells us that this set is bounded if  the vectors $\{ \psisub{k}  : 0 \le k\}$ are not restricted to any half space in $\Re^d$, for almost every initial condition  $[\varpi]$.  
The proof is postponed to the Appendix.

 \begin{proposition}
\label[proposition]{t:constGeo}
Suppose $\Prob_\varpi\{ v^\transpose \psisub{k} \ge 0 \} < 1$ for any non-zero $v\in\Re^d$.  Then, $\thetaRegion$ is compact.
\end{proposition}


\subsection{Comparison with Q-learning}

The standard Q-learning algorithm is expressed,
\begin{equation} 
\theta_{k + 1} = \theta_k + \alpha_{k + 1} \Tdiff_{k + 1}(\theta_k) \zeta_k 
\label{e:Q}
\end{equation}
where $\{\zeta_k\}$ is the sequence of $d$-dimensional eligibility vectors, typically taken as $\zeta_k = \nabla_\theta Q^\theta(Z(k))$  (which is $ \psisub{k}$ with linear function approximation),  and $\{\alpha_{k+1}\}$ the step-size sequence.   When convergent, the limit $\theta^*$ solves the so-called \textit{projected Bellman equation}  (also known as a  \textit{Galerkin relaxation}),
\begin{equation}
\Expect_\varpi\bigl[ \Tdiff_{k + 1}(\theta^*)\zeta^i_k \bigr]  =  0 \,, \qquad 1 \le i \le d \, ,
\label{e:projBE}
\end{equation}
where the expectation is in steady-state \cite{sze10,CSRL}.   


 \Cref{t:cvxQ_Galerkin} that follows shows that  	\eqref{e:ConvexQinfinity} also solves a Galerkin relaxation.  The proof follows from
\Cref{t:QCPLP}, and recognition that  in \eqref{e:ConvexQinfinityLP} we may restrict to basic feasible solutions (BFS).  

\begin{proposition}
\label[proposition]{t:cvxQ_Galerkin}
If the convex program \eqref{e:ConvexQinfinity} admits at least one optimizer,    then there is an optimizer   $\theta^*$ together with   indices $\{ i_1,\dots,i_d\} \subset \{ 1,\dots, \constD\}$ satisfying 
\begin{equation}
\Expect_\varpi\bigl[ \Tdiff_{k + 1}(\theta^*)\zeta^{i_\ell}_k \bigr]  =  0 \,, \qquad 1 \le \ell \le d \, .
\label{e:projBEcvxQ}
\end{equation}
\end{proposition}

\choreL{What does ``Possibly different" optimal solution refer to? Does it mean it's different from the solution to Galerkin relaxation in \eqref{e:projBE}?
\\
It means that not every optimizer will satisfy this.    I reworded.}

\section{Algorithms}
\label{s:algs}
In Convex Q-learning, the function $\barg: \Re^d \to \Re^\constD$ in \eqref{e:ConvexQinfinity} is replaced by its approximation via Monte Carlo
\[
\barg_N(\theta) \eqdef  \frac{1}{N} \sum_{k=0}^{N - 1}  [-\Tdiff_{k + 1}(\theta)\zeta_k].
\]


\wham{Convex Q-learning} Given the data $\{ Z(k) : 0 \le  k < N \}$ and  a pmf $\mu$ on $\zstate$, solve
\begin{equation}
\label{e:cvxQBarGTheta}
	\max_{\theta} \ \ 
	\langle  \mu, Q^\theta \rangle
	\qquad
	\st \ \ 
	\barg_N(\theta) \le 0
\end{equation}

As we increase the time horizon $N$, the variance of   the solution to \eqref{e:cvxQBarGTheta} decreases,  but the complexity of the linear program increases. 
The batch algorithms described next are designed to reduce complexity.

\subsection{Batch algorithms}
 
The two approaches below begin with the specification of intermediate times 
$T_0 = 0 < T_1 < T_2 < \cdots < T_{B-1} < T_B = N$. The parameter will be updated at these times to obtain $\{\theta_n : 0 \le n \le B\}$, initialized with $\theta_0\in\Re^d$.
Also required are two positive step-size sequences  satisfying  
\begin{equation}
	\lim_{n\to\infty}\alpha_n/\beta_n=0
	\label{e:2timescale}
\end{equation}

The empirical distribution over the $n$th batch of observations is denoted $\uppi_n$.  Hence, for any vector-valued function,  
\[
\langle \uppi_n ,  g \rangle =  \frac{1}{T_{n+1}-T_{n}}  \sum_{k=T_{n}}^{T_{n+1}-1}    g(\Phi_k )
\]
In view of \eqref{e:ConvexQinfinity},  we denote, for any $\theta\in\Re^d$,
\[
\langle \uppi_n,  \Tdiff(\theta)\zeta  \rangle \eqdef  \frac{1}{T_{n+1}-T_{n}}  \sum_{k=T_{n}}^{T_{n+1}-1}   \Tdiff_{k + 1}(\theta)\zeta_k 
\]

We  introduce a convex regularizer  $\clR_n$,    
so that the objective function at stage $n$ of the algorithm becomes 
\[
\Obj_n(\theta)  \eqdef - \langle \mu, Q^\theta \rangle  + \clR_n(\theta)
\]
and $
\Obj_n(\theta, \lambda) \eqdef  \Obj_n(\theta) -   \langle \uppi_n ,  \Tdiff(\theta )\zeta^\transpose \lambda  \rangle
$  for $\lambda\in\Re_+^{d_+}$.
Updates of $\lambda$ are obtained to approximate the Lagrange multiplier associated with   \eqref{e:ConvexQinfinity}.
Given parameter estimates $\{\theta_n :  n\ge 0 \}$, obtain a sequence of vectors in $\Re^{d_+}$ via
\[
v^{n+1} =   v^n    +  \beta_{n+1}  \bigl[  \langle \uppi_{n+1} ,  \Tdiff(\theta_{n+1})  \zeta_n   \rangle    -   v^n  \bigr]
\]
The sequence of Lagrange multiplier estimates are obtained via the recursion
$
\lambda_{n+1} = [ \lambda_{n} + \alpha_{n+1}  v^{n+1}]_+$,  
with $\lambda_0 \in\Re_+^{d_+}$ an arbitrary initial condition and $[x]_+ = \max\{0, x\}$.

\choreL{Should we find references to the two updating rule if we have space?   SPM says: We won't have space!
\\
One reference is regarding if the step size is decreasing, the solutions to the two are almost identical;
The other one is regarding there will be a difference between the two that depends on the constant learning rate.
 SPM says:   absolutely.   That is someone from the optimization literature.
\\
But which one is more preferable in practice?     SPM says:   good question!
} 
Below are two choices for parameter updates:
\wham{Batch Convex Q-learning implicit update} 
\[
\begin{aligned}
\theta_{n+1} = \argmin_\theta\Big\{  \Obj_n(\theta, \lambda_n)  + \tfrac{1}{\alpha_{n+1}}  \half  \|  \theta - \theta_n\|^2 \Big\}
\end{aligned} 
\]
It is called implicit because the solution is obtained via the fixed point equation,
\begin{equation}
\begin{aligned}
	\theta_{n+1}   =\theta_n  
	-   \alpha_{n+1} \nabla_\theta  \Obj_n(\theta, \lambda_n)\Big|_{\theta=\theta_{n+1} }
\end{aligned} 
\label{e:thetaImplicit}
\end{equation}
The explicit update is obtained by introducing a one-step delay on the right-hand side.

\wham{Batch CvxQ explicit update} 
\begin{equation}
\begin{aligned}
	&\theta_{n+1}  = \theta_n 
	-   \alpha_{n+1} \nabla_\theta \Obj_n(\theta, \lambda_n)\Big|_{\theta=\theta_{n} }
\end{aligned} 
\label{e:thetaExplicit}
\end{equation}

Assumptions on the regularizer are required to ensure convergence.   For convenience we take
\begin{equation}
\clR_n(\theta) = \kappa \{ [\langle \uppi_n  ,  \Tdiff(\theta) \zeta^\transpose \lambda  \rangle]_\_\}^2    + \epsy \| \theta \|^2
\label{e:RegN}
\end{equation}
with $\kappa, \epsy>0$ and  $[x]_- \eqdef \max(0, -x)$.   Let $\clR(\theta)$ denote its steady-state mean, and    
\[
\begin{aligned} 
\barObj(\theta, \lambda) 
\eqdef   - \langle \mu, Q^\theta \rangle 
 -   \langle \varpi ,  \Tdiff(\theta )\zeta^\transpose \lambda  \rangle + \clR(\theta)
\end{aligned} 
\]

\begin{proposition}
\label[proposition]{t:PrimalDualConvergence}
Consider either algorithm (implicit or explicit).   The algorithm is convergent to a pair $(\theta^*,\lambda^*)$   that solves the saddle point problem:
$
	\theta^* = \argmin_\theta \barObj(\theta, \lambda^*)$,
	$\lambda^* = \argmax_\lambda \min_\theta     \barObj(\theta,\lambda).
$
\end{proposition}

The proof is a standard stochastic approximation analysis,  in which the ODE approximation is precisely the   \textit{primal-dual flow}  considered in  \cite{quli19,dinjov19} and their references.

\subsection{Relative Convex Q-Learning}
\label{s:reCvxQ}

Given any pmf $\omega$ on $\zstate$, denote $H^* \eqdef Q^* - \langle \omega, Q^* \rangle$. 
Since we are subtracting a constant, we  have $\fee^*(x) = \argmin_{u\in \ustate} H^*(x, u) $.   The advantage of estimating $H^*$ is that it remains bounded for $0<\gamma<1$
  \cite{devmey22}.  The extension of the preceding theory to this \textit{relative Q-function} is straightforward, beginning with

\begin{proposition} 
	\label[proposition]{t:relpcharact}
	For any positive pmfs  $\mu$ and $\omega$ on $\zstate$ and  positive scalar $\delta > 0$,  $H^*$ solves the  convex program,
	\begin{equation}
		\begin{aligned}
			\max_H &\ \ 
			\langle \mu,  H \rangle
			\\
			\st &\ \ 
			H(z) \le c(z) + \gamma P_u \uH(x) - \delta \langle \omega, H\rangle \,, 	
			\\		
		    & \qquad \textit{for every} \  z = (x, u)\in \zstate.
		\end{aligned}
		\label{e:RelaCvxQ}
	\end{equation}
\end{proposition}
 This motivates one version of relative convex Q-learning:

\wham{Relative CvxQ} Given the data $\{ Z(k) : 0 \le k < N \}$ and  probability measures $\mu$ and $\omega$ on $\zstate$, solve
\begin{equation}
	\label{e:cvxRQTheta}
		\max_{\theta} \ \ \langle  \mu, H^\theta \rangle
		\qquad
		\st \ \ \barg_N(\theta) \le 0 
\end{equation}
where $\barg_N$ is defined as in \eqref{e:cvxQBarGTheta}, with $\Tdiff_{k + 1}$ replaced by the relative temporal difference:
\begin{align*}
	\TdiffR_{k+1}(\theta) &\eqdef -H^\theta(Z(k)) 
	+  c(Z(k))  + \gamma\uH^\theta(X(k+1))
	\\
	&\qquad- \delta \langle \omega, H^\theta\rangle.
\end{align*}
Formulation of batch algorithms is also straightforward.

\section{Rates of Convergence}
\label{s:vartheory}

We consider here the rate of convergence for the basic algorithm \eqref{e:cvxQBarGTheta},  whose solution is denoted $\theta_N$.   Subject to mild  conditions we establish that $\theta_N \to \theta^*$ with probability one as $N\to\infty$,  where $\theta^*$ solves \eqref{e:ConvexQinfinity};  one assumption is that the solution is unique.     It is more challenging to establish bounds on the mean-square error (MSE) 
$\Expect[\|\tiltheta_N\|^2]$,  with $\tiltheta_N = \theta_N -\theta^*$.      In fact, we have not been able to find a deterministic $N_0$ for which   $\Expect[\|\tiltheta_N\|^2] <\infty$ for $N\ge N_0$. 
We fix $r>0$ satisfying $|\theta^*_i| < r$ for each $i$,   and let $\theta_N^r$ denote the $L_\infty$ projection of the solution of \eqref{e:cvxQBarGTheta}  to the region $\thetaRegion_r =\{\theta\in\Re^d :    |\theta_i| \le r  \,, \ 1\le i\le d\}$.     We set $\theta_N^r = 0$ if  the convex program   \eqref{e:cvxQBarGTheta}  is unbounded or infeasible.

While the sequence $\{\theta_N^r\}$ cannot be represented as the output of any recursive algorithm,  key results from 
stochastic approximation theory would suggest that the MSE should decay as $O(1/N)$.
We verify that this rate of convergence holds, and obtain finer results:
\begin{equation}
	\begin{aligned}
		N  \Expect[  \tiltheta^r_N    ( \tiltheta^r_N ) ^\transpose ]   &=   \SigmaTheta + O\bigl(\tfrac{1}{\sqrt{N} } \bigr)
		\\
		\sqrt{N} \tiltheta_N^r & \darrow N(0, \SigmaTheta) \,, \qquad N\to\infty
	\end{aligned} 
	\label{e:MSECLT}
\end{equation}
where the convergence is in distribution in the second limit,  and $\tiltheta^r_N = \theta^r_N -\theta^*$.  The matrix $\SigmaTheta \ge 0$ is identified in  \Cref{t:cvxQcovariance}
after a few preliminary results.

For each $\theta\in\Re^d$ and $\fee\in\Upphi$ denote
\[
\barg_N(\theta,\fee) \eqdef   -  \frac{1}{N} \sum_{k=0}^{N - 1}  g_k(\theta, \fee)  \, ,
\] 
where $g_k(\theta, \fee) =  -\Tdiff_{k+1}(\theta, \fee)\zeta_k$  (recall \eqref{e:Tdiff_theta+fee}).

\begin{proposition}
	\label[proposition]{t:QCPLPMCE}
	Any solution to the convex program \eqref{e:cvxQBarGTheta} is also a solution to the linear program
	\begin{equation}
		\label{e:cvxQThetaLP}
		\begin{aligned}
			\max_\theta &\ \ 
			\langle  \mu, Q^\theta \rangle
			\\
			\st &\ \ 
			\barg_N^i(\theta, \fee)  \le 0 \,, \ \ 1 \le i \le \constD, \ \ \fee \in \Upphi \, .
		\end{aligned}
	\end{equation}
	Moreover, a vector 
	$\theta\in\Re^d$ is a BFS if and only if the following two properties hold: 
	1.\  there is a set $ \clI_+^\theta = \{ j_1,\dots, j_d\} \subset \{1,\dots, \constD\}$  such that 
	\[
	\begin{aligned}
		\barg_N^i(\theta, \fee^\theta)   & = 0 \,, \ \   i\in  \clI_+^\theta
		\\
		&  <0   \,, \ \   i\not\in  \clI_+^\theta
	\end{aligned} 
	\]
	where $\fee^\theta$ is the $Q^\theta$-greedy policy,  and 2.\ There is no other $Q^\theta$-greedy policy:
	 if $\fee'\in\Upphi$ satisfies
	\[
	\fee'(x) \in \argmin_u  Q^\theta(x,u)\,,\qquad \textit{for all $x\in\state$,}
	\]
	then $\fee'(x) = \fee^\theta(x)$ for all $x$.    
\end{proposition}

 \wham{Proof.}
The proof of the LP characterization \eqref{e:cvxQThetaLP}   is identical to the proof of \Cref{t:QCPLP}.    The characterization of a BFS is a consequence of the proof:   By definition, if $\theta$ is a BFS then there are exactly $d$ pairs $\{ \fee^i, j_i : 1\le i\le d \}$ for which the constraints are tight,  meaning $\barg_N^{j_i} (\theta, \fee^i)    = 0$.   We also have $\barg_N^{j_i} (\theta, \fee^\theta) \ge \barg_N^{j_i} (\theta, \fee^i)    = 0$, so that by feasibility we must also have      $\barg_N^{j_i} (\theta, \fee^\theta) =0$ for each $i$.      Since there are exactly $d$ active constraints, we must have $\fee^\theta =  \fee^i$ for each $i$.  
 \qed

\Cref{t:QCPLPMCE} provides the ingredients required to establish convergence.

It is assumed that $\theta^*$ is unique, from which it follows that it is also a BFS.    We let $\clI_+ = \{ j_1,\dots, j_d\} $ denote the set of $d$ indices for which 
$\barg_N^i (\theta, \fee^\theta) =0$  when $\theta=\theta^*$ and $i\in\clI_+$.

\choreL{
	I see why the reviewer are confused with the `+' sign. The definition of  $\barA^+$ is not mentioned.
	According to the appendix, $\barA$ is defined as the steady state mean, while $\barA^+$ is obtained by selecting $d$ rows form $\barA$ and the indices are from $\clI_+$.
	\\
	Sean, maybe you have better idea on how to introduce these terms?
	\\
	I simply forgot the "+" when I defined $A_k^+$!  Fixed.
	}

A $d$-dimensional stochastic process $\barW_N$ is constructed as an average,
\[
\barW_N = \frac{1}{N} \sum_{k=1}^N  W_k
\]
in which the $d$-dimensional stochastic process $\{W_k \}$ is obtained in the following steps.    First,  construct a $d\times d$ matrix $A_k^+$, whose   $i,j$ element is given by 
\[
[A_k^+]_{i,j} \eqdef 
\bigl[ -\psisub{k-1}
+ \gamma  \upsisub{k}^{\theta^*} \bigr]_j \zeta_{k-1}^{j_i}    
\]
and let
$\barA^+=   \Expect_\varpi\bigl[A^+_k]$, where the expectation is in steady-state.    Let $\beta_k^+ $ be the $d$-dimensional vector whose $i$th component is equal to  $ \csub{k-1}  \zeta_{k-1}^{j_i}  $,   whose steady state mean is the $d$-dimensional vector  $\barbeta^+  =  \Expect_\varpi[\beta_k^+ ]$.    
We then take 
\[ 
W_k  =  \beta_k^+ - \barbeta^+   - [ A_k^+ - \barA^+]\theta^*
\]
This is analogous to the ``disturbance sequence'' that arises in variance analysis of standard Q-learning algorithms \cite{devmey17b,CSRL}.   

	 
\begin{proposition}
\label[proposition]{t:cvxQcovariance}
Suppose that  the linear program \eqref{e:cvxQThetaLP} has a unique optimizer $\theta^*$.  
Then, $\displaystyle \lim_{N\to\infty} \theta^r_N =\lim_{N\to\infty} \theta_N = \theta^*$ with probability one, and \eqref{e:MSECLT} holds.      
	
The covariance matrix may be expressed  
\[
\begin{aligned}
		\SigmaTheta & =  [\barA^+]^{-1} \Sigma_W   ([\barA^+]^{-1})^\transpose  \,,
		\\
		\textit{with} \ \         \Sigma_W & =  \lim_{N\to\infty} N  \Expect[  \barW_N (\barW_N)^\transpose ]   
	\end{aligned} 
	\]
\end{proposition}

  Invertibility of $\barA^+$  is a consequence of the fact that $\theta^*$ is a BFS. 

The form of the covariance $\SigmaTheta$ is identical in form to the minimal covariance identified in the averaging theory of Polyak and Ruppert for stochastic approximation  (see historical notes in \cite[Ch.~8]{CSRL}).   It is likely that it is also minimal here (the batch algorithms  may satisfy a CLT, but the covariance matrix will dominate  $\SigmaTheta$).

\choreL{$\barA^+$ is also mentioned in this paragraph.
\\
Sorry, that was a typo}

\wham{Proof of \Cref{t:cvxQcovariance} (main concepts).}
The full proof follows from the more general  \Cref{t:randLP} found in the Appendix.
  We present here the main ideas through  a heuristic, which is conveniently explained for the more general setting of  convex constraints\footnote{Extension of variance theory to general convex constraints  may be valuable when relaxing the assumption that $\ustate$ is finite, or when convex regularizers are introduced.}.

Consider the minimization of a linear objective $\objv^\transpose \theta $ subject to convex constraints $\barg_N(\theta)\le 0$,   and let $\theta^*$ denote the minimizer for $N=\infty$: 
as in \Cref{t:cvxQcovariance} there is a limiting function $\barg$,  and   $\theta^*$ minimizes $\objv^\transpose \theta $ subject  $\barg(\theta)\le 0$.
The heuristic is based on the convex program
\begin{align}
	\label{e:varQP}
\theta_N^* = \argmin \objv^\transpose \theta   \quad \text{subject to }  \ \  \barg(\theta)  \le  b_N^*
\end{align}
with $b_N^* \eqdef  \barg(\theta^*) - \barg_N(\theta^*) $.   The remainder of the proof is based on the fact that $\Expect[\| b_N^*\|^2] = O(1/N)$,  along with sensitivity theory for convex programming.

In analysis of 	\eqref{e:varQP} the definition of the set $\clI_+$ remains the same,     characterized in terms of  the dual under mild conditions.   
Let $\lambda^*\in\Re^{\constD}$ denote the Lagrange multiplier associated with the optimizer $\theta^*$. 
If sensitivity is strictly positive for active constraints then $\clI_+=\{ i : \lambda^*_i>0\}$.
\qed

\begin{figure}[h]
	\begin{center}
		\centerline{\includegraphics[width=.75\columnwidth]{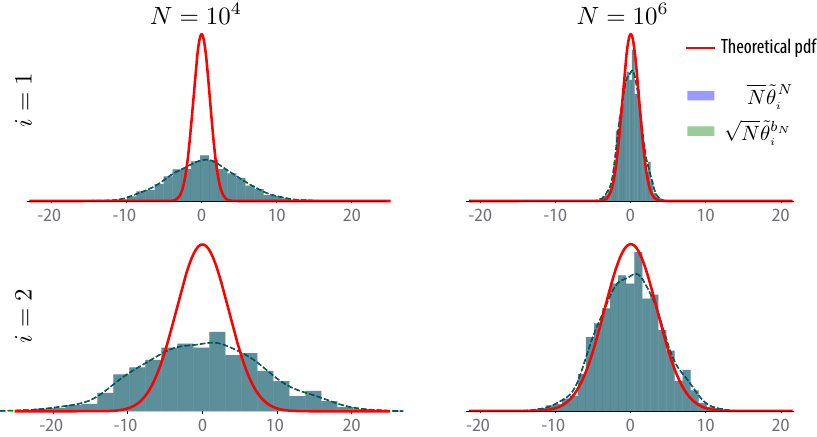}}
		\caption{Histogram plots of $\sqrt{N}\tiltheta_N$ and $\sqrt{N}\tiltheta^*_N$. }
		\label[figure]{f:varTheory}
	\end{center}
 \vspace{-.4cm}
\end{figure}

 \wham{Example}
 
To illustrate the CLT in a simple setting,    consider a quadratic program on $\Re^2$ with 10 constraints: $\theta\in\Re^2$, $\barg: \Re^2 \to \Re^{10}$, and
$
\barg_N(\theta) =N^{-1} \sum_{k = 0}^{N - 1} g(\theta + \Delta_k) $,   with $\{\Delta_k\}$ i.i.d.,  $ N(0, I)$.    
The function $g:\Re^2 \to \Re^{10}$ is the quadratic
$
g(\vartheta) = (a^\transpose \vartheta)\odot (a^\transpose \vartheta) + b^\transpose \vartheta - \One$, 
where $a, b \in \Re^{2\times 10}$, $\One$ is a vector of ones, and $\odot$ denotes element-wise multiplication.

Three quadratic programs were constructed by replacing the constraint in \eqref{e:varQP} with $\barg_N(\theta) \le 0$, $\barg(\theta) \le 0$, and $\barg_N(\theta) \le b^*_N$. The solutions to each of these quadratic programs are denoted   $\theta_N$, $\theta^*$, and $\theta_N^*$ respectively.

The values of $a$, $b$, and $\objv$ were obtained by sampling independently from a normal distribution.   The results that follow show typical results for one set of values.

The dual variable $\lambda^*$ associated with  $\theta^*$ was also obtained. Exactly  two of ten inequality constraints were found to be tight, and exactly two values of $\lambda^*$ were strictly positive.

The theoretical Gaussian density was compared to histograms obtained from   repeated independent experiments.  In each of 100 runs,  the corresponding errors were recorded: 
 $\sqrt{N}\tiltheta^*_N \eqdef \sqrt{N} [\theta^*_N - \theta^*]$ and $\sqrt{N}\tiltheta_{N} \eqdef \sqrt{N}[\theta_{N} - \theta^*]$ .   
 Results shown in \Cref{f:varTheory} are consistent with the CLT.  Approximations are good even for   $N = 10^4$, and nearly perfect for $N\ge 10^6$.

\section{Application to Inventory Control}
\label{s:numerics}

We survey here results from experiments on a classical inventory control problem focusing on two topics: 
(i) Stability and consistency of CvxQ and relative CvxQ; 
(ii) Comparison of convergence rates with Q-learning and CvxQ.\footnote{The experiments go beyond the theory because the state space is not finite.  We believe that numerical results are valuable for testing the boundaries of the theory.   We know from experience that some  disagree, so we respectfully ask the reviewers to accept our preferences.}

\begin{subequations}

\wham{Preliminaries}
Consider the  inventory model with inventory level $X(k)\in\state=\Re$ (a negative value indicating backlog),   depletion rate $\beta>0$,  and stocking decision $U(k)\in  \ustate = \{ 0,  1\}$.  The MDP model has cost defined by parameters $ c^+,c^- > 0$ and evolution equation,
\begin{equation}
\begin{aligned}
&X(k+1) = X(k) -  [\beta+ W(k+1)  ]+ U(k) 
		\\
&		c(x,u) = \max(c^+ x, -c^-  x)   \,, \  \ x\in\Re\, ;
	\end{aligned}
\label{e:inventoryModel}
\end{equation}
  $\{W(k)\}$ is  i.i.d.\  with zero mean and  finite variance $\sigma^2_W$.

\end{subequations}



An optimal policy is of the threshold form:  
\begin{align}
	\label{e:infee}
	\fee (x; \barr) = \ind\{ x \le -\barr \}.
\end{align}
For small   $\beta>0$  
the optimal threshold is approximated by $
\barr^\dagger \eqdef \log(1 +c^+/c^-)/\varrho 
$
with
$\varrho$ the positive solution to 
$
	\sigma^2_W \varrho^2 / 2 - \beta\varrho - \gamma = 0
$.  See \cite[Ch.~7]{CTCN} for background.

\begin{figure}[h]
	\begin{center}
		\centerline{\includegraphics[width=.75\columnwidth]{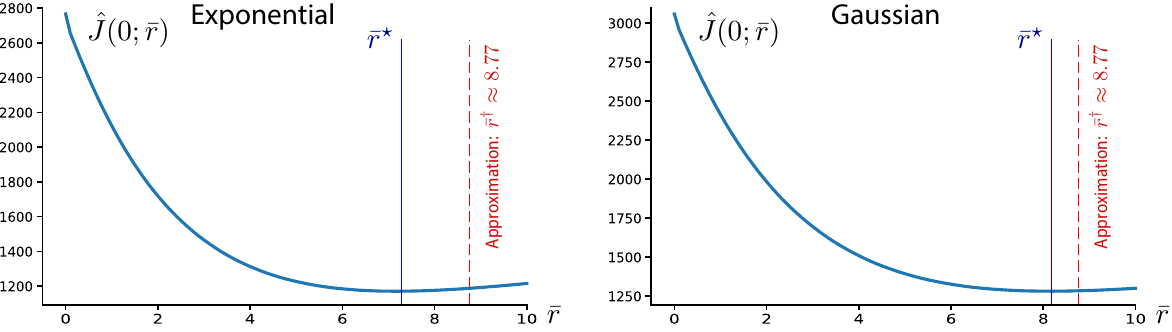}}
		\caption{Numerical Estimation of the Optimal Inventory Level.}
		\label{fig:numEstBarR}
	\end{center}
	\vspace{-0.4cm}
\end{figure}

We set $\gamma = 0.99$, $\beta = 0.1$, $c^- = 1$, $c^+ = 10$ and $\sigma^2_W=1$ in the experiments that follow, giving $\barr^\dagger \approx 8.77$.



\wham{Numerical Estimation of $\bfmath{\bar r^*}$.}

To evaluate the results obtained using CvxQ we estimated the optimal threshold via Monte-Carlo,  evaluating a range of 100 values of $\barr$ evenly spaced in the interval $[0, 10]$.       
We used common randomness for the entire range of thresholds, meaning that we constructed a single i.i.d.\ family  
$\{ W_k^{i}: 1 \le i \le \barN,   1 \le k \le N  \}$,    with   $\barN = 2\times 10^4$;
 in each experiment we chose  initial condition $X^i(0) = 0$.   
 For each threshold $\barr$ we obtain the sample path $\{X^i(k;\barr):  0\le k<N\}$,  $ N = 10^4$, and from this   an estimate of the discounted total cost by truncating   and averaging:   
\[
\haJ(x;\barr) = \frac{1}{\barN} \sum_{i = 1}^{\barN} \sum_{k = 0}^{N - 1} \gamma^k c(X^i(k;\barr))   \,, \ \ \barr \in [0, 10]   \, .
\]

The plots in \Cref{fig:numEstBarR} show results for two choices of disturbance,  $W_k \sim N(0,1)$ and    $W_k \sim \text{Exp}(1) -1$.    In each case, the value of  $\haJ(0;\barr) $  at $\barr^\dagger$ closely matches the empirical minimum $\barr^*$.   


\begin{figure*}
	\begin{center}
		\centerline{\includegraphics[width=1\hsize]{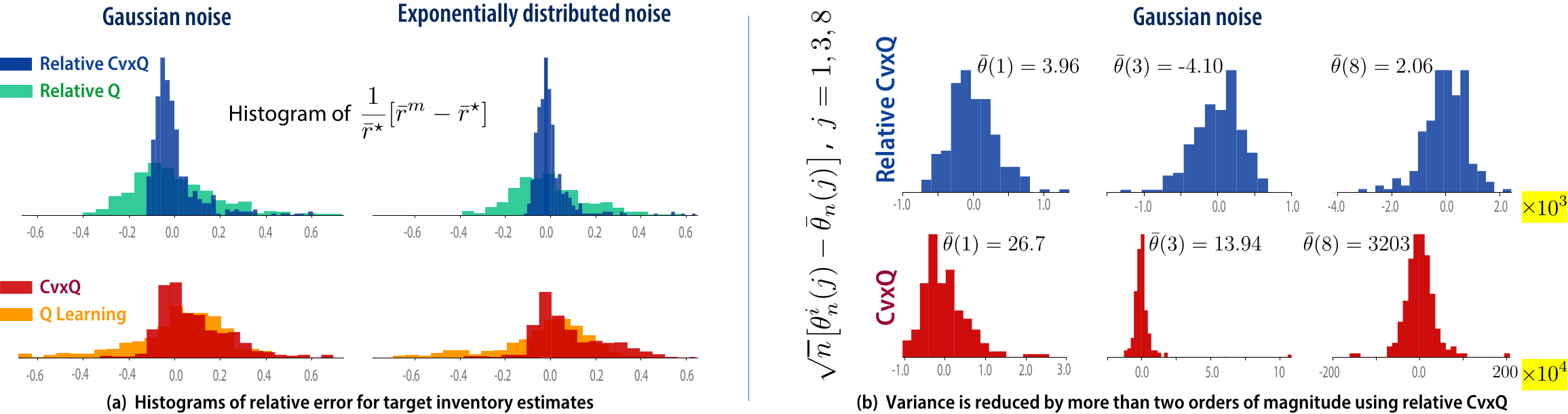}}
		\caption{Performance Comparison between CvxQ, CvxRQ Q-learning and relative Q-learning.}
		\label{fig:cvxQRQInvenLevelEst}
	\end{center}
	\vspace{-0.4cm}
\end{figure*}

%
%



\subsection{Algorithm Comparisons}

  CvxQ, relative CvxQ, Q-learning, and relative Q-learning  were applied to approximate the optimal policy under the same modeling assumptions used to obtain the plots in  \Cref{fig:numEstBarR}.

 \wham{Details of Implementation}

Theory surveyed in \cite[Ch.~7]{CTCN} tells us that the optimal value function  $J^*(x)\eqdef\min_u Q^*(x,u)$ is convex,  and is approximated by $c(x)/ (1-\gamma)$ for large $x$.    
This motivated the choice of  8-dimensional  linear function class  $\{Q^\theta = \theta^\transpose \psi: \theta \in \Re^d\}$  using   the continuously differentiable functions
		$\psi(z) = [\psi^\prime(x)\ind\{ u =0 \} ; \psi^\prime(x)\ind\{ u =1 \}]^\transpose$ where
$
\psi^\prime(x) =  [\xi_1(x);\xi_2(x); x;1]$,  with $\xi_i(x) = 
 (|x| + (e^{-\delta_i |x|} - 1))\ind\{ x \ge 0 \} / \delta_i$.
      The values   $\delta_1 = 0.5$, $\delta_2 = 0.1$, gave good results.   
 

The following input sequence was used for training:
\begin{equation*}
	\begin{aligned}
		\fee(u | x) &\eqdef P(U(k ) = u | X(k) = x)
		\\
		&= \epsy P_E(u) + (1 - \epsy)\ind\{ u = \fee(x) \}
	\end{aligned}
\end{equation*}
with  $\epsy = 0.9$,  $P_E$   uniform   on $\ustate$,   and   $\fee$   a threshold policy that gave good performance. 

We chose $\{\zeta_k\}$ in  \eqref{e:cvxQBarGTheta}, \eqref{e:cvxRQTheta} 
to be the indicator functions:
\begin{equation*}
	\zeta^i_k = \ind \{x^{i} \le X(k) \le x^{i+1} \} \,, \ \ 1 \le i \le \constD
\end{equation*}
We found   that choosing $x^i$ to be evenly spaced in the range $[-28, 28]$ gave  good performance.  
The results that follow applied this approach using $\constD = 200$.




Several algorithms were tested:    Q-learning using the recursion \eqref{e:Q} with $\zeta_k = \psisub{k}$ and constant step size $\alpha_{k + 1} = 10^{-3}$;   
relative Q-learning algorithm,  defined   by replacing $\Tdiff_{k + 1}(\theta_k)$ with $\TdiffR_{k + 1}(\theta_k)$ in \eqref{e:Q};   
CvxQ \eqref{e:cvxQBarGTheta} and relative CvxQ \eqref{e:cvxRQTheta} using the convex solver from the Python library  CVXPY.  
 
Independent experiments were conducted with common disturbance  $\{W(k)\}$ in each run (using the model \eqref{e:inventoryModel}).   
 In each of  $M = 100$ independent runs with time horizon $N=10^4$,  the following data was collected:     the final parameter estimate $\theta^m$,  and the estimate $\barr^m$ of the optimal threshold, defined as the minimal solution
 $Q^{\theta^m}(x, 0) \le  Q^{\theta^m}(x, 1) $ for all $x\ge \barr^m$.      Selected results are illustrated in \Cref{fig:cvxQRQInvenLevelEst}.
Shown on the left hand side are histograms of the relative error $\{ [\barr^m -\barr^*]/\barr^* : 1\le m\le M \}$.      We see that relative CvxQ offers the lowest variance in these experiments.  

The right hand side of \Cref{fig:cvxQRQInvenLevelEst}  shows histograms of components of the scaled parameter error,   $\sqrt{n}[\theta^*_m(i) - \bartheta(i)]$  for several  $i$,  with $\bartheta$ the average of $\theta^*_m$.  Only results for   the CvxQ algorithms are shown.   
 The  CLT appears to hold, with relative CvxQ giving much lower variance.


\section{Conclusions}

This paper lays out new approaches to reinforcement learning design and analysis.   Important gaps include the augmentation of CvxQ to impose bounds on performance of resulting policies, techniques for analysis in non-stationary settings (motivated by the need for more efficient exploration),   and development of algorithms for more advanced function approximation architectures (e.g., neural networks and RKHS).    
It is important to improve the numerical performance of   CvxQ,   build bridges with  logistic Q-learning and recent approaches, such as the interacting particle approach of \cite{jostagmehmey22},   and find ways to make comparisons in terms of both training efficiency and policy performance.

\clearpage

\normalem
\bibliographystyle{abbrv}

\def\cprime{$'$}\def\cprime{$'$}



\clearpage

\begin{appendices}
	
	\renewcommand{\thetheorem}{A.\arabic{theorem}}

	\section{Technical Proofs}


\wham{Proof of \Cref{t:constGeo}}  The proof is similar to the characterization of boundedness in the deterministic setting \cite{lumehmeyneu22}.   We highlight the steps and minor differences here.  

If the constraint region is unbounded, then for each $m \ge 1$, there exists $\theta^m$ such that $\| \theta^m \| \ge m$, and
\[
	\frac{1}{\|\theta^m\|} \Expect\bigl[-\Tdiff_{k + 1}(\theta^m) | \clF_k\bigr] \le 0
\]
Let $\hatheta$ denote a   limit point of the sequence  $\{\hatheta^m \eqdef \theta^m / \|\theta^m\| : m\ge 1\}$:   $\hatheta =  \lim_{i \to \infty} \hatheta^{m_i}$ for a subsequence $\{m_i\}$.
By the definition of $\uQ^{\theta^m}$  and elementary arguments we obtain   
	\begin{equation*}
		\begin{aligned}
		0 &\ge 
			\lim_{i \to \infty} \Expect[-\widehat\Tdiff_{k + 1}(\hatheta^{m_i})| \clF_k] 
			\\
			& =
			\Expect  \bigl[Q^{\hatheta}(Z(k))  - \gamma\uQ^{\hatheta}(X(k + 1)) | \clF_k\bigr]  			\\
			& \ge
			\Expect \bigl[Q^{\hatheta}(Z(k))  - \gamma Q^{\hatheta}(Z(k + 1)) | \clF_k\bigr]  
		\end{aligned}
	\end{equation*}

Non-positivity of the right hand side may be expressed,  
\[
	\Expect \bigl[\hatheta^\transpose (\psi_{(k)}  - \gamma \psi_{(k+1)}) | \clF_k\bigr]  = -\epsy_k
\]
where $\{\epsy_k : k \ge 1\}$ is a non-negative sequence.      On multiplying each side by   $\gamma^k$, and following the steps in \cite{lumehmeyneu22},     
\[
		\begin{aligned}
			\Expect \bigl[\hatheta^\transpose (\gamma^{k + i} \psi_{(k+i)}  &- \gamma^{k +i+1} \psi_{(k+i+1)}) | \clF_{k}\bigr]
			\\
			&= -\gamma^{k+i}\Expect[\epsy_{k+i} | \clF_k]
		\end{aligned}
\]
With $k$ fixed,  sum from $i = 0$ to $N\ge 1$ to obtain
	\begin{align*}
	-\gamma^{k} \sum_{i = 0}^N\gamma^{i}\Expect[\epsy_{k+i} | \clF_k] 
& =
		\sum_{i = 0}^N \Expect \bigl[\hatheta^\transpose (\gamma^{k + i} \psi_{(k+i)}  - \gamma^{k +i+1} \psi_{(k+i+1)}) | \clF_{k}\bigr]
		\\
		&=\gamma^{k} \hatheta^\transpose  \psisub{k}  -   \gamma^{k}\Expect \bigl[  \hatheta^\transpose \gamma^{N + 1} \psisub{N+1} | \clF_{k}\bigr]
	\end{align*}
This gives the desired result on letting $N\uparrow\infty$:
	\begin{align*}
		\gamma^{k}\hatheta^\transpose  \psisub{k}   = -\gamma^{k}\sum_{i = 0}^\infty\gamma^{i}\Expect[\epsy_{k+i} | \clF_k] \le 0
	\end{align*}
	\qed
	
\wham{Proof of \Cref{t:cvxQcovariance}.}
To ease notation in analysis of the basic Convex Q-learning algorithm we pose an abstraction. 
Suppose that $\bfmY = \{ Y_k : k\ge 0\}$ is a Markov chain on a finite state space $\ystate$ that is uni-chain and aperiodic.  In application to \Cref{t:cvxQcovariance} we have $Y_k = (X_k; U_k; X_{k+1})$.

The assumptions on $\bfmY$ imply that   there is a unique invariant pmf $\varpi$ such that $\lim_{k\to\infty} \Prob\{Y_k = y' \mid Y_0 = y\} = \varpi(y')$ for all $y,y'\in\ystate$, where the rate of convergence is necessarily geometric.    For integers $1< d< \constD$ and functions $A\colon \ystate \to \Re^{\constD\times d}$ and $\beta\colon \ystate \to \Re^{\constD }$ we denote $A_k = A(Y_k)$,  $\beta_k = \beta(Y_k)$,  and the averages
$
\barA_N = N^{-1}\sum_{k=1}^N A_k$,  $  
\barbeta_N = N^{-1}\sum_{k=1}^N \beta_k
$.
In view of \Cref{t:QCPLPMCE} the convex program \eqref{e:cvxQTheta} can be expressed as a linear program:
\begin{equation}
	\theta_N =  \min\{   \objv^\transpose \theta \,,
	\quad \st \ 
	\barA_N\theta  \le \barbeta_N  \} \, .    
	\label{e:randLP}
\end{equation}

Convergenceof  $\{ \theta_N \}$ in an almost sure sense requires uniqueness of the solution to the limiting LP,
\begin{equation}
	\theta^* \eqdef \objv^\transpose \theta \,,
	\qquad \st \ \ 
	\barA\theta  \le \barbeta \, .    
	\label{e:randLPinfinity}
\end{equation}
with $\barA = \sum\varpi(y) A(y)$,   $\barbeta = \sum\varpi(y) \beta(y)$.

For mean-square convergence we require projection, as assumed in \Cref{t:QCPLPMCE}:   fix $r>0$ satisfying $|\theta^*_i| < r$ for each $i$,   and let $\theta_N^r$ denote the $L_\infty$ projection of the solution of \eqref{e:randLP}  to the region $\thetaRegion_r =\{\theta\in\Re^d :    |\theta_i| \le r  \,, \ 1\le i\le d\}$.     We set $\theta_N^r = 0$ if  the LP   \eqref{e:randLP}  is unbounded or infeasible.

The \textit{fundamental theorem of linear programming} tells us that uniqueness of $\theta^*$ implies that it is a BFS.   
Let $ \clI_+$ denote the set of $d$ indices for  which the constraint is tight:
\begin{equation}
	[\barA\theta^*  - \barbeta ]_\ell = 0    \,,   \ \     \ell \in  \clI_+ = \{ i_1,\dots, i_d\} \, ,
	\label{e:BFS}
\end{equation}
and $[\barA\theta^*  - \barbeta ]_i <0 $  otherwise.
Letting $\barA^+$  denote the $d\times d$ matrix  whose $\ell$th row is equal to $i_\ell$,
 and  $\barbeta^+ = (\barbeta_{i_1}; \dots; \barbeta_{i_d})\in\Re^d$, 
we have $\theta^* =  [\barA^+ ]^{-1} \barbeta^+$.

  The  $d\times d$ matrix $\tilA_k^+  $ is obtained by extracting the same rows from  $\tilA_k \eqdef A_k -\barA$,   and the $d$-dimensional vector $\tilbeta_k^+$ is obtained by extracting these elements from $  \tilbeta_k = \beta_k -\barbeta$.

 \choreL{$i_j$? Do we mean $j_i \in\clI_+$? How about $j_i\in \clI_+$ \\
 It was a mess.  OK now}

A CLT   for the error $\theta_N^r -\theta^*$  is provided in \Cref{t:randLP}.  Its proof is based on  justifying the approximation $\theta_N^r \approx \theta^* +  [\barA^+ ]^{-1}     \barW_N $ with 
\[
\barW_N  \eqdef  \frac{1}{N} \sum_{k=1}^N   W_k\,,  \qquad W_k \eqdef  \tilbeta_k^+ -  \tilA_k^+ \theta^*    \, .
\]


\begin{proposition}
	\label[proposition]{t:randLP}
If the solution to \eqref{e:randLPinfinity} is unique, then $\displaystyle\lim_{N\to\infty} \theta^r_N =\lim_{N\to\infty} \theta_N = \theta^*$ with probability one.
	
The mean-square error admits the  approximation,
\[
	N  \Expect[   ( \theta^r_N -\theta^*)   ( \theta^r_N -\theta^*) ^\transpose ]   =   \SigmaTheta + O\bigl(\tfrac{1}{\sqrt{N} } \bigr)
	\]
	where $\SigmaTheta =  [\barA^+]^{-1} \Sigma_W   ([\barA^+]^{-1})^\transpose    $,  with
	\[
	\Sigma_W  =  \lim_{N\to\infty} N  \Expect[  \barW_N (\barW_N)^\transpose ]   
	\]
	Moreover,  $\sqrt{N}  ( \theta^r_N -\theta^*)  \darrow N(0,\SigmaTheta)$, where the convergence is in distribution.
\end{proposition}

The proof is based on a  version of the CLT from  \cite{ste01a}:
\begin{proposition}
	\label[proposition]{t:ste01a}
	For any function $g\colon\ystate\to\Re$,  let $\tilg_k = g(Y_k) -\varpi(g)$  (centered to have zero steady-state mean),  and let $Z_N =  N^{-1/2} \sum_{k=1}^N \tilg_k$.     
	Then, for any continuous function $F\colon\Re\to\Re$ with polynomial growth,  and each initial condition $Y_0$,
	\begin{equation}
		\lim_{N\to\infty}  \Expect [ F(Z_N)]   = \Expect [ F(Z)]   \,,  \quad  Z\sim N(0, \gamma^2_g)
		\label{e:ste01a}
	\end{equation}
	where $\gamma^2_g \ge 0$ is the asymptotic variance of $\{ \tilg_k \}$  (the variance appearing in the CLT for $\{Z_N\}$).   If $\gamma^2_g > 0$ then the continuity assumption on $F$ may be relaxed.  
\end{proposition}

\wham{Proof of \Cref{t:randLP}}  
The proof will follow from  the bound
\begin{equation}
	\Expect[\| \clE_N^\uptheta\|^2] = O(1/N^2)\,, \quad
	\clE_N^\uptheta    \eqdef
	\theta^r_N -\theta^*  - [\barA^+]^{-1} \barW_N 
	\label{e:randLPkeyBdd}
\end{equation}

Subject to \eqref{e:BFS} the following conclusions hold when  $(\barA_N,\barbeta_N)$ is sufficiently close to $(\barA,\barbeta)$:     
$\theta_N$ is also unique, with tight constraints specified by $ \clI_+$, and   
\begin{equation}
	\theta_N = [ \barA_N^+]^{-1} \barbeta_N^+
	\label{e:thetaN_BFS}
\end{equation}   

Let   $\clS_N \in\clF$ denote the event 
\[
\clS_N  =  \{  \| \barA_N -\barA \|_F < \delta \quad \textit{and } \quad \| \barbeta_N  -\barbeta\| < \delta
\}
\]
where $\| \varble \|_F$ denotes the Frobenious norm,   and  
$(\Omega,\clF)$ is the underlying probability space.    It is assumed that $\delta>0$ is chosen sufficiently  small so that   \eqref{e:thetaN_BFS} holds 
on this event, and   moreover $\| \clE_N\|_F\le 1/2$ with $\clE_N = I - \barA_N^+[ \barA^+]^{-1}   $.     This final bound is imposed since we will apply the following:
\begin{equation}
	\begin{aligned}
		\barA^+&[\barA_N^+]^{-1}   = [I - \clE_N]^{-1}  =  I+ \clE_N   +  [I -\clE_N]^{-1} \clE_N^2
		\\
		&   \textit{with} \ \    \|  [I -\clE_N]^{-1}   \|_F  \le 2  \quad    \textit{ on the event $\clS_N $}  
		\\
		&   \textit{and} \ \    \Expect[ \| \clE_N^2 \|_F^2] = O(1/N^2)
	\end{aligned} 
	\label{e:barAinvN}
\end{equation}
where the $L_2$ bound follows from \Cref{t:ste01a}.
The following representation is also useful:
\begin{equation}
	\barW_N = \barbeta_N^+ -\barbeta^+  +   \clE_N  \barbeta^+    
	\label{e:WclE}
\end{equation}

Applying \eqref{e:ste01a} and Markov's inequality gives for $p\ge 1$,
\begin{equation}
	\begin{aligned}
		\Prob\{\clS_N^c\}   &
		\le   \Prob \{  N^{p/2} \| \barA_N -\barA \|_F^p \ge  N^{p/2}\delta^p \}  
		\\
		&\qquad +    \Prob \{ N^{p/2}\| \barbeta_N  -\barbeta\|^p \ge N^{p/2} \delta^p\}
		\\
		& \le  N^{- p/2} \delta^{-p}    B_N
	\end{aligned} 
	\label{e:ste01aProb}
\end{equation}
The bounded sequence  $\{B_N : N\ge 1\}$    is a linear combination of expectations of the form $ \Expect [ F(Z_N)] $ as appearing in 
\eqref{e:ste01a};    \Cref{t:ste01a} is applied for several functions,  $g(y) = A_{i,j}(y)$ and $\beta_i(y)$ for every possible $i,j$,  and
$F(y) =   |y|^p$.   It follows that  $\Prob\{\clS_N^c\} $  vanishes quickly as $N\to\infty$.

Write
$
\clE_N^\uptheta    =   \clE_N^0   \ind_{\clS_N}  +   \clE_N^\uptheta  \ind_{\clS_N^c}   
$ (recall  \eqref{e:randLPkeyBdd}).  
The second moment of the second term decays to zero faster than $O(1/N^q)$ for any integer $q\ge 2$:  This is a consequence of \eqref{e:ste01aProb}, the bound   $\|  \theta_N^r -\theta^*  \|_\infty \le 2r$ (with $r$ used in the projection),  and Chebyshev's inequality:
\[
\Expect[ \|  [\barA^+]^{-1} \barW_N  \| \ind_{\clS_N^c}   ]^2   \le \Expect[ \|  [\barA^+]^{-1} \barW_N  \|^2]  \Prob\{\clS_N^c\} 
\]

Of course, $ \clE_N^0 = \clE_N^\uptheta   $.  The change of notation is because, on the event  $\clS_N$,  we may apply   \eqref{e:thetaN_BFS} to obtain 
\[
\clE_N^0   = [ \barA^+]^{-1}  \bigl[  \barA^+ [ \barA_N^+]^{-1} \barbeta_N^+     -   \barbeta^+  -   \barW_N   \bigr]
\]
where we have from \eqref{e:barAinvN},
\[
\begin{aligned}
	& \barA^+[\barA_N^+]^{-1}  -  (I+ \clE_N )  =    Q_N\,,  \ \ Q_N \eqdef
	[I -\clE_N]^{-1} \clE_N^2
	\\
	&  \qquad    \|Q_N\|_F  \le  2  \|\clE_N^2 \|_F  \quad \text{(on the event $\clS_N$).}
\end{aligned} 
\]
Applying \eqref{e:WclE} then gives
\[
\begin{aligned} 
	\barA^+ [ \barA_N^+]^{-1} \barbeta_N^+        -   \barbeta^+   
	&  =      (1+ \clE_N) \barbeta_N^+      + Q_N \barbeta_N^+  -\barbeta^+
	\\
	&= \barW_N +   \clE_N (\barbeta_N^+ -\barbeta^+)    + Q_N\barbeta_N^+
\end{aligned} 
\]
and on applying \Cref{t:ste01a},
\[
\begin{aligned}
	\Expect[ \| & \clE_N^0\|^2   \ind_{\clS_N} ]
	\\
	&\le  \kappa \Bigl(\Expect[ \|  \clE_N  (\barbeta_N^+ -\barbeta^+)  \|^2]  +  
	\Expect[ \|  Q_N \|_F^2 \| \barbeta_N^+\| ^2 \ind_{\clS_N} ]  \Bigr)
	\\
	& \le O(1/N^2)  \,,  \qquad \textit{with $\kappa = 2 \|  [ \barA^+]^{-1}\|_F $.}
\end{aligned} 
\]
This establishes the bound in \eqref{e:randLPkeyBdd}
and completes the proof.   
\qed	

 \end{appendices}

\end{document}